\documentclass{article}
\PassOptionsToPackage{numbers,compress}{natbib}
 \usepackage[sglblindworkshop, final]{neurips_2025}
\workshoptitle{MLxOR: Mathematical Foundations and Operational Integration of Machine Learning for Uncertainty-Aware Decision-Making}

\usepackage[utf8]{inputenc} 
\usepackage[T1]{fontenc}    
\usepackage{hyperref}       
\usepackage{url}            
\usepackage{booktabs}       
\usepackage{amsfonts}       
\usepackage{nicefrac}       
\usepackage{microtype}      
\usepackage{xcolor}         

\title{Differentiable Optimization for \\Deep Learning-Enhanced DC Approximation \\of AC Optimal Power Flow}

\author{%
  Andrew Rosemberg  \\   NSF AI Institute for Advances in Optimization\\ Georgia Institute of Technology\\ Atlanta, GA \\
  \texttt{arosemberg3@gatech.edu}
  \\ \And
  Michael Klamkin \\ NSF AI Institute for Advances in Optimization\\ Georgia Institute of Technology\\ Atlanta, GA \\
  \texttt{klam@gatech.edu}
  \\ \And
  Pascal Van Hentenryck \\ NSF AI Institute for Advances in Optimization\\ Georgia Institute of Technology\\ Atlanta, GA \\
  \texttt{pvh@gatech.edu}
}

\usepackage{graphicx}
\usepackage{amsmath}
\usepackage{amssymb}
\usepackage{booktabs}
\usepackage{bm}

\newcommand{\gff}{g^{\text{ff}}}
\newcommand{\gft}{g^{\text{ft}}}
\newcommand{\gtf}{g^{\text{tf}}}
\newcommand{\gtt}{g^{\text{tt}}}
\newcommand{\bff}{b^{\text{ff}}}
\newcommand{\bft}{b^{\text{ft}}}
\newcommand{\btf}{b^{\text{tf}}}
\newcommand{\btt}{b^{\text{tt}}}

\newcommand{\pgmin}{\underline{\mathbf{p}}^{\text{g}}}
\newcommand{\pgmax}{\overline{\mathbf{p}}^{\text{g}}}
\newcommand{\qgmin}{\underline{\mathbf{q}}^{\text{g}}}
\newcommand{\qgmax}{\overline{\mathbf{q}}^{\text{g}}}

\newcommand{\dvamin}{\underline{\Delta} \theta}
\newcommand{\dvamax}{\overline{\Delta} \theta}

\newcommand{\EDGES}{\mathcal{E}} 

\newcommand{\LOADS}{\mathcal{L}}
\newcommand{\GENERATORS}{\mathcal{G}}
\newcommand{\NODES}{\mathcal{N}}
\newcommand{\PG}{\mathbf{p}^{\text{g}}} 
\newcommand{\QG}{\mathbf{q}^{\text{g}}} 
\newcommand{\PD}{\mathbf{p}^{\text{d}}} 
\newcommand{\QD}{\mathbf{q}^{\text{d}}} 
\newcommand{\GS}{{g}^{\text{s}}} 
\newcommand{\BS}{{b}^{\text{s}}} 
\newcommand{\VM}{\mathbf{v}} 
\newcommand{\VA}{\bm{\theta}} 
\newcommand{\PF}{\mathbf{p}^{\text{f}}} 
\newcommand{\QF}{\mathbf{q}^{\text{f}}} 
\newcommand{\PT}{\mathbf{p}^{\text{t}}} 
\newcommand{\QT}{\mathbf{q}^{\text{t}}} 

\newcommand{\lambdaP}{\lambda^{\text{p}}}

\newcommand{\lambdaPf}{{\lambda}^{\text{pf}}}

\newcommand{\muPgMin}{\underline{\mu}^{\text{pg}}}
\newcommand{\muPgMax}{\bar{\mu}^{\text{pg}}}

\newcommand{\muPfMin}{\underline{\mu}^{\text{pf}}}
\newcommand{\muPfMax}{\bar{\mu}^{\text{pf}}}

\newcommand{\muAngleDiffMin}{\underline{\mu}^{\theta}}
\newcommand{\muAngleDiffMax}{\bar{\mu}^{\theta}}

\usepackage{xcolor,stackengine}
\newcommand\radbox[1]{%
  \fboxsep=2pt
  \fboxrule=.75pt
  \def\tmp{\displaystyle\strut #1}
  \def\shadow{\makebox[0pt]{$\tmp$}}
  \stackengine{0pt}{%
    \stackengine{0pt}{%
      \textcolor{red}{\fbox{$\phantom{\tmp}$}}%
    }{\color{white}}{O}{c}{F}{F}{L}%
  }{$\tmp$}{O}{c}{F}{F}{L}
}

\usepackage{subcaption}
\usepackage{wrapfig}
\usepackage{float}
\newfloat{model}{thp}{lop}\floatname{model}{Model}
\begin{document}

\maketitle

\begin{abstract}
The growing scale of power systems and the increasing uncertainty introduced by renewable energy sources necessitates novel optimization techniques that are significantly faster and more accurate than existing methods. The AC Optimal Power Flow (AC-OPF) problem, a core component of power grid optimization, is often approximated using linearized DC Optimal Power Flow (DC-OPF) models for computational tractability, albeit at the cost of suboptimal and inefficient decisions. To address these limitations, we propose a novel deep learning-based framework for network equivalency that enhances DC-OPF to more closely mimic the behavior of AC-OPF. The approach utilizes recent advances in differentiable optimization, incorporating a neural network trained to predict adjusted nodal shunt conductances and branch susceptances in order to account for nonlinear power flow behavior. The model can be trained end-to-end using modern deep learning frameworks by leveraging the implicit function theorem. Results demonstrate the framework's ability to significantly improve prediction accuracy.
\end{abstract}

\section{Introduction}
Power systems optimization has gained significant attention in recent years, driven by growing emphasis on expansion, the integration of renewable energy sources, and the general need for cleaner, more sustainable energy solutions \cite{s2022, pozo2012three, pozo2012chance,barry2022risk}.
A fundamental component of optimal power grid control is the constrained optimization problem known as AC Optimal Power Flow (AC-OPF, Model \ref{model:acopf}, see Table \ref{tab:notations} for notations), which includes the nonconvex AC Power Flow (AC-PF) equations \eqref{model:acopf:kirchhoff:active}$-$\eqref{model:acopf:ohm:reactive:to} to model how power is transmitted through the grid's transmission lines and transformers, engineering constraints such as transmission line thermal limits \eqref{model:acopf:thrmbound:from}$-$\eqref{model:acopf:thrmbound:to}, and minimizes total power generation cost \eqref{model:acopf:obj}. Obtaining accurate OPF solutions is critical to many downstream tasks across time scales, including real-time risk-aware market clearing \cite{tam2011real,e2elr},
  day-ahead security-constrained unit commitment ({SCUC})~\cite{sun2017decomposition},
  transmission switching optimization~\cite{OTS},
  and expansion planning~\cite{verma2016transmission}.  The main challenge precluding more widespread use of AC-OPF in practice lies in the non-linear and non-convex nature of the physics and engineering constraints, making solving time scale poorly as network size and uncertainty grows. Consequently, AC-(O)PF is often not directly utilized in practice \cite{o2011recent}. Instead, transmission system operators (TSOs) opt for simplified problems such as DC-(O)PF \cite{o2011recent, eldridge2016marginal}, where assumptions are made to approximate the physics using only linear functions, facilitating the use of more computationally efficient convex optimization techniques. While this approach is tractable, the approximations currently in use lead to sub-optimal decisions that incur significant operational inefficiency and thus increased cost \cite{rosemberg2021assessing}.

\begin{wrapfigure}{r}{0.56\textwidth}
    \centering
    \includegraphics[width=0.98\linewidth]{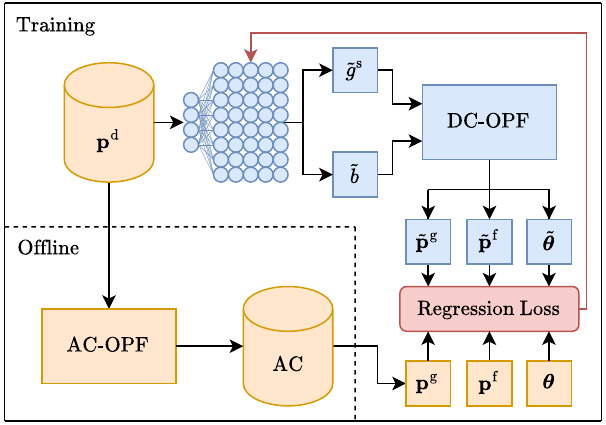}
    \caption{The proposed framework for deep learning-based DC-OPF adjustment.}
    \label{fig:diagram}
\end{wrapfigure}

Recent research has focused on leveraging deep learning to address these challenges. Most approaches (as comprehensively reviewed in \cite{litreview}) employ ``proxy models" — surrogate models trained via supervised learning, often with additional loss function terms to reduce constraint violations. These models directly approximate the power demand to optimal solution mapping. While such models have shown promising results, achieving high-fidelity predictions with minimal constraint violations for AC-OPF \cite{van2021machine, chen2022learning, piloto2024canos,klamkin2024bucketized}, they typically lack guarantees of optimality or feasibility, particularly when applied to scenarios outside their training distribution. This limitation poses a significant challenge in ensuring reliable operations in the face of uncertainty.

This work addresses this challenge by proposing a novel deep learning framework that leverages recent advances in the differentiation of constrained optimization with respect to problem parameters \cite{agrawal2019differentiable} in order to \textit{learn how to conditionally modify the problem data of DC-OPF such that its optimal solution mimics locally the behavior of AC-OPF}. Figure \ref{fig:diagram} summarizes the approach. First, a dataset of power demands $\PD$ is sampled and the corresponding AC-OPF solutions are obtained using a non-linear programming solver such as Ipopt \cite{ipopt}. {The neural network model is created by appending a differentiable DC-OPF layer to a feed-forward neural network that predicts model parameters. The feed-forward portion of the neural network takes as input the power demand $\mathbf{p}^\text{\normalfont d}$ and outputs two vectors: adjusted nodal shunt conductances $\tilde{g}^\text{\normalfont s}$ and adjusted branch susceptances $\tilde{b}$, which are then used to formulate an ``adjusted'' DC-OPF. The final layer solves this adjusted DC-OPF to obtain its optimal solution $\tilde{\mathbf{p}}^\text{\normalfont g}$, $\tilde{\mathbf{p}}^\text{\normalfont f}$, $\tilde{\VA}$} which is the overall model's prediction. Adjusting $\GS$ and $b$ allows the model to account for power losses and variations in effective line impedance that would otherwise be neglected in a standard DC-OPF formulation (which assumes loss-less power flow and linearized Kirchhoff laws), making it possible for the adjusted DC-OPF to better approximate the nonlinear characteristics of AC-OPF while retaining interpretability.

\section{Related Work}
The need for accurate yet computationally tractable models for power system optimization is well-documented. Particularly, \citet{rosemberg2021assessing} underscores the economic consequences of relying on simplified models, examining the trade-offs associated with convex relaxations and approximations in hydrothermal dispatch planning. Notwithstanding, the Federal Energy Regulatory Commission (FERC) report \cite{o2011recent} details current practices in the U.S., where transmission system operators often adopt linear relaxations, despite their potential to yield infeasible or inefficient solutions.

Classical approaches to power system reduction trace back to \cite{ward1949equivalent}, which introduced network reduction models that remain foundational in simplifying power systems. However, these methods lack the flexibility needed to address the variability and uncertainty inherent in modern, renewable-integrated grids.
Recent works have proposed new methods to overcome these limitations. For example, \cite{ribeiro2023equivalent} explores network reduction models using nonlinear basis functions, leveraging the fact that basis function fitting can be formulated as linear constraints, offering a computationally efficient alternative. Additionally, \cite{litreview} provides a comprehensive review of deep learning applications in OPF, highlighting the growing interest in data-driven approaches that aim to enhance the scalability and robustness of OPF solutions. \citet{constante2024ac} proposes an approach that first solves the bilevel problem of finding the best DC-OPF parameters to exactly match AC-OPF solutions for each instance in a dataset, then fitting a neural network to map AC-OPF loads to these parameters. While this approach minimizes the error between the AC and DC formulations, the requirement for exact bilevel solutions severely limits scalability.

\section{Technical Approach}
\paragraph{Problem Setup}
The learning problem can be posed as finding the neural network weights $\omega$ that approximately solve the bilevel optimization problem:
                    \begin{equation*} 
                        \operatorname*{argmin}_{\omega} \;\mathbb{E}_{\PD} \left[
                        \begin{aligned}
                            &\operatorname{MSELoss}\big(\begin{bmatrix}
                            {\mathbf{p}}^\text{g} \\ {\mathbf{p}}^\text{f} \\ {\VA}
                        \end{bmatrix},\;\begin{bmatrix}
                            \PG_\text{AC}\\ \PF_\text{AC} \\ \VA_\text{AC}
                        \end{bmatrix}\big)
                    \quad\text{where} \; 
                        \begin{aligned}[t]
                            \begin{bmatrix}
                            {\mathbf{p}}^\text{g} \\ {\mathbf{p}}^\text{f} \\ {\VA}\end{bmatrix} &\in \operatorname{argmin}\; \operatorname{DC-OPF}(\PD;\; {g}^\text{s},\, {b}) \\
                        \end{aligned} \\[-1.5em]
                        & \quad\quad\quad\quad\quad\quad\quad\quad\quad\quad\quad\quad\quad\quad\quad\quad\quad\quad\quad\quad\enspace\;\text{\rm s.t.} \quad\enspace\GS,\, b = \operatorname{NN}_\omega(\PD)
                    \end{aligned}
                    \right]
                \end{equation*}

In order to find the mapping from power demand to optimally adjusted nodal shunt conductances and adjusted branch susceptances, $\mathbf{p}^\text{\normalfont d} \rightarrow (\tilde{g}^\text{\normalfont s}, \tilde{b})$, one needs to solve this complicated bilevel problem, i.e. taking into account the optimal solution of the lower level problem (DC-OPF) for any possible choice of upper level solution. The optimality of the lower level can be characterized using the Karush-Kuhn-Tucker (KKT) conditions which consist of guaranteeing: (a) Primal Feasibility, (b) Stationarity, (c) Dual Feasibility, and (d) Complementary Slackness. The exact mathematical forms of these conditions for DC-OPF is included in \ref{sec:dckkt}.

\paragraph{Offline AC-OPF Data Generation}
The proposed method is a ``semi-self-supervised'' learning approach; it relies on a dataset of power demands and corresponding AC-OPF solutions, but solves DC-OPF in-the-loop. Thus, before training begins, a number of AC-OPF instances must be solved offline.
These will be used in an ordinary regression loss, minimizing the distance between the adjusted DC-OPF optimal solution and the AC-OPF optimal solution.

\paragraph{Differentiable Optimization Layers}
The core of the proposed method is the use of a differentiable optimization layer to enable backpropagation through the solving of the DC-OPF, then differentiating the solution with respect to problem parameters which were predicted by the feed-forward neural network portion. Although automatic differentiation can in principle be used to differentiate through the iterations of an LP solver, this technique scales poorly with the number of iterations, resulting in massive computational graphs. It is also cumbersome to implement in practice since most modern solvers do not natively integrate with automatic differentiation systems. Instead of differentiating through the solution iterates, the {differentiable optimization layer approach uses the implicit function theorem to solve for the gradient of the solution map with respect to problem data in terms of the optimal solution}, which is obtained in the forward pass by calling a solver \cite{agrawal2019differentiable,diffopt}.

Specifically, the sensitivities are calculated by applying the implicit function theorem applied to the KKT conditions of the optimization problem. These KKT conditions encapsulate the solution of the DC-OPF problem as a system of equations, i.e.
$F(\mathbf{x}, \phi) = 0$,
where
\(\mathbf{x} = (\PG, \PF, \VA, \lambda)\) represents the decision variables of the OPF, including power generation (\(\PG\)), power flow (\(\PF\)), voltage angles (\(\VA\)), and Lagrange multipliers (\(\lambda\)); and \(\phi = (\tilde{g}^\text{\normalfont s}, \tilde{b})\) represents the adjustable problem data: nodal shunt conductances and branch susceptances. The exact KKT conditions are given in Section \ref{sec:dckkt}.

This capability allows for the proposed framework to focus on learning how to account for the nonlinearities only rather than also trying to learn how to solve the DC-OPF itself (an end-to-end proxy has to learn both simultaneously). This is achieved by having the neural network predict DC-OPF problem data (nodal shunt conductance and branch susceptance) then using a linear optimization solver to obtain the corresponding DC-OPF solution -- rather than predicting the solution directly.

\paragraph{Baseline Methods}
The proposed method is compared to the standard proxy approach \cite{litreview} where a feed-forward neural network is trained using MSE loss to learn the map from power demand to AC-OPF solutions directly. The proxy and the proposed method are also compared to using the DC-OPF solution directly.
The methods are evaluated using solution accuracy; that is, the $\text{L}_1$ distance between the predicted $\tilde{\mathbf{p}}^\text{g}$, $\tilde{\mathbf{p}}^\text{f}$, $\tilde{\VA}$ and the true (i.e. found by Ipopt \cite{ipopt}) $\PG$, $\PF$, $\VA$.

\section{Results}
This section includes experimental results comparing the baseline optimization proxy method, the current approach of using the DC-OPF solution itself, and the proposed neural-adjusted DC-OPF (hereafter referred to as DC2AC). Additional details on experiment settings are given in Section \ref{sec:expdetail}

\subsection{Accuracy Analysis}

\begin{figure}[t]
  \centering

  \begin{subfigure}{0.32\linewidth}
    \centering
    \includegraphics[width=\linewidth]{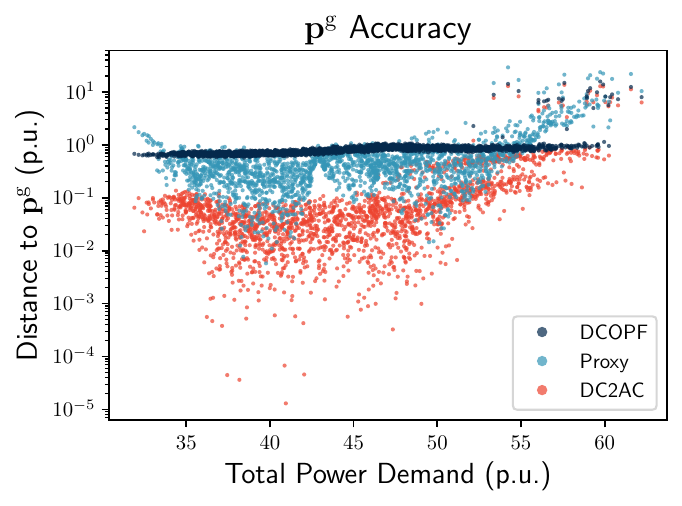}
    \caption{$\|\tilde{\mathbf{p}}^{\text{g}}-\mathbf{p}^{\text{g}}\|_1$}
    \label{fig:dcac_pg}
  \end{subfigure}\hfill
  \begin{subfigure}{0.32\linewidth}
    \centering
    \includegraphics[width=\linewidth]{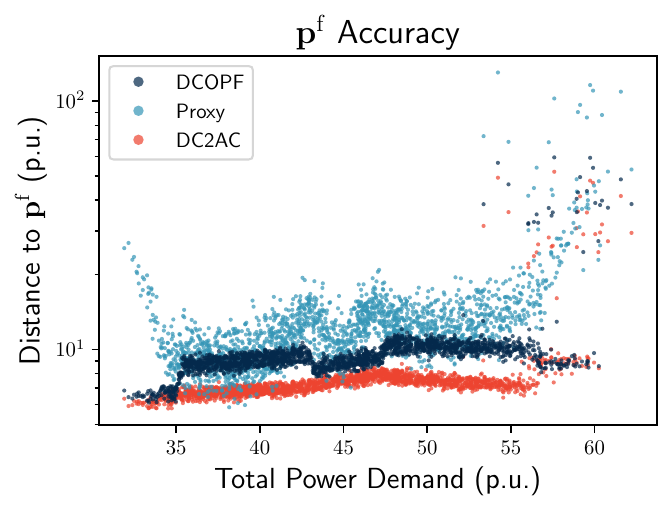}
    \caption{$\|\tilde{\mathbf{p}}^{\text{f}}-\mathbf{p}^{\text{f}}\|_1$}
    \label{fig:dcac_pf}
  \end{subfigure}\hfill
  \begin{subfigure}{0.32\linewidth}
    \centering
    \includegraphics[width=\linewidth]{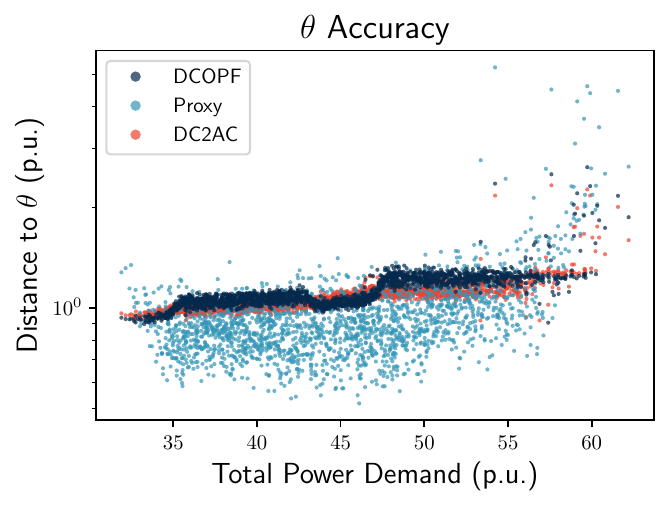}
    \caption{$\|\tilde{\VA}-\VA\|_1$}
    \label{fig:dcac_va}
  \end{subfigure}

  \caption{Accuracy vs. total demand across methods.}
  \label{fig:accuracy_panels}
\end{figure}

Table \ref{tab:accuracy} summarizes the validation set accuracy of each approach, broken down by type of decision variable. Evidently, the DC2AC approach outperforms both baselines for $\PG$ and $\PF$, though the proxy is the best for $\VA$. The columns of Table \ref{tab:accuracy} are further visualized in Figure \ref{fig:accuracy_panels}.

Notably, across all decision variables, all methods struggle in the high-load region corresponding to the power system being congested. Interestingly, the DC-OPF often outperforms the proxy in the high-load region, and DC2AC typically improves slightly on DC-OPF. For $\PG$ and $\PF$, especially in the lower-load portion of the domain, DC2AC is usually the best method of the three, achieving the lowest errors compared to the true AC solution. Overall, DC2AC outperforms all other methods on $95\%$ of the validation set samples for $\PG$ and $97\%$ of samples for $\PF$, but only $13\%$ of samples for $\VA$. Remarkably, DC2AC is able to almost exactly reproduce the AC-optimal $\PG$ in mid-range cases, achieving errors as low as $10^{-5}$. The relative flexibility of the proxy is shown to help to learn the more complicated $\VA$ variables, with DC2AC struggling to significantly outperform the vanilla DC-OPF. We hypothesize that this is due to the fact that the DC-OPF formulation we use neglects reactive power \cite[Appendix A.2.3]{opfgenerator}, a key limitation that is, evidently, hard to overcome for DC2AC.

Figure \ref{fig:conv} shows the convergence plots for the proxy and DC2AC learning approaches, where solid lines denote validation loss and dashed lines denote training loss. Note that the networks have identical initialization, and that the DC2AC training was stopped early due to lack of improvement.
The figure shows that DC2AC converges extremely quickly, mostly within the first and second epochs.

\section{Conclusion}
This paper introduced DC2AC, an interpretable deep learning approach to approximate AC-OPF solutions. The core idea is to use a differentiable optimization layer to solve a parametrized DC-OPF in the forward pass, leveraging the implicit function theorem to enable end-to-end training.
Overall, the results show that the DC2AC approach is promising, with clear directions for further improvement. We hypothesize that since output of DC2AC is directly the (parametrized) DC-OPF solution, it struggles to capture the active-reactive power relationship, explaining the poor performance on the $\VA$ variables. Future work includes exploring
the use of different approximations of AC-OPF to address this and integration with GPU solvers \cite{gpulp} to speed up training. 

\clearpage
\section*{Acknowledgements}
  This material is based upon work supported by the National Science Foundation Graduate Research Fellowship under Grant No. DGE-2039655 and the National Science Foundation Artificial Intelligence Institute for Advances in Optimization under Grant No. 2112533. Any opinions, findings, and conclusions or recommendations expressed in this material are those of the author(s) and do not necessarily reflect the views of the National Science Foundation. \nocite{PACE}
{
\bibliographystyle{unsrtnat}
\bibliography{egbib}
}

\clearpage
\appendix
\section{Formulations and Notations}

This section includes the detailed mathematical programming formulations referenced throughout the paper. Table \ref{tab:notations} summarizes the notations, Model \ref{model:dcopf} presents the DC-OPF, \ref{model:acopf} presents the AC-OPF, and \ref{model:d-dcopf} presents the dual of DC-OPF. The KKT conditions of DC-OPF are also included in section \ref{sec:dckkt}.
\begin{table}[!ht]
\centering
\caption{OPF Data Notations}
    \begin{tabular}{c|c|c}
        \toprule
        Description & Size & Symbol \\
        \midrule
        Set of buses & $-$ & $\NODES$ \\
        Set of edges & $-$ & $\EDGES$ \\
        Set of loads & $-$ & $\LOADS$ \\
        Set of generators & $-$ & $\GENERATORS$ \\
        Reference bus voltage angle & 1 & $\VA_\text{ref}$ \\
        Active power load & $|\LOADS|$ & $\PD$ \\
        Reactive power load & $|\LOADS|$ & $\QD$ \\
        Active power generation & $|\GENERATORS|$ & $\PG$ \\
        Reactive power generation & $|\GENERATORS|$ & $\PG$ \\
        Active power flow from & $|\EDGES|$ & $\PF$ \\
        Active power flow to & $|\EDGES|$ & $\PT$ \\
        Reactive power flow from & $|\EDGES|$ & $\QF$ \\
        Reactive power flow to & $|\EDGES|$ & $\QT$ \\
        Voltage magnitude & $|\NODES|$ & $\VM$ \\
        Voltage angle & $|\NODES|$ & $\VA$ \\
        Nodal load shedding & $|\NODES|$ & $\varphi$ \\
        Indices of branches leaving each bus & $|\NODES|$ & $\EDGES_i$ \\
        Indices of branches entering each bus & $|\NODES|$ & $\EDGES^R_i$ \\
        Indices of generators at each bus & $|\NODES|$ & $\GENERATORS_i$ \\
        Indices of loads at each bus & $|\NODES|$ & $\LOADS_i$ \\
        Nodal shunt conductance & $|\NODES|$ & $\GS$ \\
        Nodal shunt susceptance & $|\NODES|$ & $\BS$ \\
        Voltage magnitude lower bound & $|\NODES|$ & $\underline{\VM}$ \\
        Voltage magnitude upper bound & $|\NODES|$ & $\overline{\VM}$ \\
        Minimum voltage angle difference & $|\EDGES|$ & $\dvamin$ \\
        Maximum voltage angle difference & $|\EDGES|$ & $\dvamax$ \\
        Branch thermal limit & $|\EDGES|$ & $\overline{S}$ \\
        Minimum active power generation & $|\GENERATORS|$ & $\underline{\PG}$ \\
        Maximum active power generation & $|\GENERATORS|$ & $\overline{\PG}$ \\
        Minimum reactive power generation & $|\GENERATORS|$ & $\underline{\QG}$ \\
        Maximum reactive power generation & $|\GENERATORS|$ & $\overline{\QG}$ \\
        Generator cost coefficient & $|\GENERATORS|$ & $c$ \\
        Load shedding cost & $1$ & $c_\varphi$ \\
        From bus index for each branch & $|\EDGES|$ & $i$ \\
        To bus index for each branch & $|\EDGES|$ & $j$ \\
        Branch conductance & $|\EDGES|$ & $g$ \\
        Branch susceptance & $|\EDGES|$ & $b$ \\
        From-side branch conductance & $|\EDGES|$ & $\gff$ \\
        From-to branch conductance & $|\EDGES|$ & $\gft$ \\
        To-from branch conductance & $|\EDGES|$ & $\gtf$ \\
        To-side branch conductance & $|\EDGES|$ & $\gtt$ \\
        From-side branch susceptance & $|\EDGES|$ & $\bff$ \\
        From-to branch susceptance & $|\EDGES|$ & $\bft$ \\
        To-from branch susceptance & $|\EDGES|$ & $\btf$ \\
        To-side branch susceptance & $|\EDGES|$ & $\btt$ \\
        \bottomrule
    \end{tabular}
\label{tab:notations}
\end{table}

\begin{model}[!ht]
  \caption{DC Optimal Power Flow (DC-OPF)}
  \label{model:dcopf}
  \begin{subequations}
      \begin{align*}
          \min_{\PG, \PF, \VA, \bm\varphi} \quad&
              \sum_{i \in \NODES} c_\varphi \varphi_i+\sum_{j \in \GENERATORS_{i}} c_j \PG_j  \\
          \text{s.t.} \quad
          & \sum_{j\in\GENERATORS_i}\PG_j - \sum_{e \in \mathcal{E}_{i}}  \PF_{e} + \sum_{e \in \mathcal{E}^R_{i}} \PF_{e} + \varphi_i
            = \sum_{j\in\LOADS_i}{\PD_j} + \radbox{\GS_i}
              & \forall i &\in \NODES  & [\lambdaP]\\
          & \PF_{e} = {-}\radbox{b_{e}}(\VA_{i} - \VA_{j})
              & \forall e = (i, j) &\in \EDGES & [\lambdaPf] \\
        & \dvamin_{e} \leq \VA_i - \VA_j \leq \dvamax_{e}
              & \forall e = (i, j) &\in \EDGES
               &[\muAngleDiffMin,\muAngleDiffMax]\\
          & \VA_\text{ref} = 0  \\
          & \underline{\PG_i} \leq \PG_i \leq \overline{\PG_i}
              & \forall i &\in \GENERATORS
               & [\muPgMin,\muPgMax] \\
          & {-}\overline{S_{e}} \leq  \PF_{e} \leq \overline{S_{e}}
              & \forall e &\in \EDGES
               & [\muPfMin,\muPfMax]
      \end{align*}
  \end{subequations}
\end{model}

\begin{model}[!ht]
    \caption{Dual of DC-OPF}
    \label{model:d-dcopf}
    \begin{subequations}
        \begin{align}
            \max_{\substack{\lambda^\text{p},\lambda^\text{pf},\underline{\mu}^\theta,\overline{\mu}^\theta,\\\underline{\mu}^\text{pg},\overline{\mu}^\text{pg},\underline{\mu}^\text{pf},\overline{\mu}^\text{pf}}} \quad
            & \sum_{i \in \NODES} \lambdaP_{i} (\radbox{\GS_{i}} + \varphi_i + \sum_{j \in \LOADS_{i}}  \PD_{j} )
            + \sum_{i \in \NODES} \sum_{j \in \GENERATORS_{i}} \left(
                \pgmin_{j} \muPgMin_{j} - \pgmax_{j} \muPgMax_{j}
            \right) \notag\\
            & + \sum_{e \in \mathcal{E}} \left(
                  \dvamin_{e} \muAngleDiffMin_{e}
                - \dvamax_{e} \muAngleDiffMax_{e}
                - \bar{s}_{e} \muPfMin_{e}
                - \bar{s}_{e} \muPfMax_{e}
            \right)
            \label{model:d-dcopf:obj} \\
            \text{s.t.} \quad
            & \lambdaP_{i} + \muPgMin_{g} - \muPgMax_{g} = c_{g}
                & \forall i \in \NODES, \forall g &\in \GENERATORS_{i}
                \label{model:d-dcopf:pg}\\
            &  -\lambdaP_{i} + \lambdaP_{j} - \lambdaPf_{e} + \muPfMin_{e} - \muPfMax_{e}
                = 0
                & \forall e=(i, j) &\in \EDGES
                \label{model:d-dcopf:pf}\\
            &  \sum_{e \in \EDGES_{i}^{+}} \left(\muAngleDiffMin_{e} -  \radbox{b_{e}} \lambdaPf_{e} \right)
                + \sum_{e \in \EDGES^{-}_{i}} \left(\radbox{b_{e}} \lambdaPf_{e} - \muAngleDiffMax_{e}\right)
                = 0
                & \forall i &\in \NODES
                \label{model:d-dcopf:va}\\
            & \muAngleDiffMin, \muPgMin, \muPfMin \geq 0\\
            & \muAngleDiffMax, \muPgMax, \muPfMax \geq 0
        \end{align}
    \end{subequations}
\end{model}

\begin{model}[!ht]
  \caption{AC Optimal Power Flow (AC-OPF)}
  \label{model:acopf}
  \begin{subequations}
      \begin{align}
          \min_{\substack{\PG, \QG, \VM, \VA\\\PF, \QF, \PT, \QT }} \quad
            & \sum_{i \in \NODES} \sum_{j \in \GENERATORS_{i}} c_j \PG_j \label{model:acopf:obj} \\
          \text{s.t.} \quad 
          & \sum_{j\in\GENERATORS_i}\PG_j - \sum_{j\in\LOADS_i}\PD_j - \GS_i \VM_i^2 = \sum_{e \in \mathcal{E}_{i}}  \PF_{e} + \sum_{e \in \mathcal{E}^R_{i}} \PT_{e}
              & \forall i &\in \NODES \label{model:acopf:kirchhoff:active} \\
          & \sum_{j\in\GENERATORS_i}\QG_j - \sum_{j\in\LOADS_i}\QD_j + \BS_i \VM_i^2 = \sum_{e \in \mathcal{E}_{i}} \QF_{e} + \sum_{e \in \mathcal{E}^R_{i}} \QT_{e}
              & \forall i &\in \NODES \label{model:acopf:kirchhoff:reactive} \\
          & \PF_{e} = \gff_{e}\VM_i^2 + \gft_{e} \VM_i \VM_j \cos(\VA_i-\VA_j) + \bft_{e} \VM_i \VM_j \sin(\VA_i-\VA_j)
              & \forall e = (i,j) &\in \EDGES \label{model:acopf:ohm:active:from} \\
          & \QF_{e} = -\bff_{e} \VM_i^2 - \bft_{e}\VM_i \VM_j \cos(\VA_i-\VA_j) + \gft_{e} \VM_i \VM_j \sin(\VA_i-\VA_j)
              & \forall e = (i,j) &\in \EDGES \label{model:acopf:ohm:reactive:from} \\
          & \PT_{e} = \gtt_{e}\VM_j^2 + \gtf_{e} \VM_i \VM_j \cos(\VA_i-\VA_j) - \btf_{e} \VM_i \VM_j \sin(\VA_i-\VA_j)
              & \forall e = (i,j) &\in \EDGES\label{model:acopf:ohm:active:to} \\
          & \QT_{e} = -\btt_{e} \VM_j^2 - \btf_{e}\VM_i \VM_j \cos(\VA_i-\VA_j) - \gtf_{e} \VM_i \VM_j \sin(\VA_i-\VA_j)
              & \forall e = (i,j) &\in \EDGES \label{model:acopf:ohm:reactive:to} \\
          & (\PF_{e})^2 + (\QF_{e})^2 \leq \overline{S_{e}}^2
              & \forall e &\in \EDGES
              \label{model:acopf:thrmbound:from} \\
          & (\PT_{e})^2 + (\QT_{e})^2 \leq \overline{S_{e}}^2
              & \forall e &\in \EDGES
              \label{model:acopf:thrmbound:to} \\
          & \dvamin_{e} \leq \VA_i - \VA_j \leq \dvamax_{e}
              & \forall e = (i,j) &\in \EDGES
              \label{model:acopf:angledifference} \\
          & \VA_\text{ref} = 0 \label{model:acopf:slackbus} \\
          & \underline{\VM_i} \leq \VM_i \leq \overline{\VM_i}
              & \forall i &\in \NODES
              \label{model:acopf:vmbound} \\
          & \pgmin_{i} \leq \PG_i \leq \pgmax_{i},\;\qgmin_{i} \leq \QG_i \leq \qgmax_{i}
              & \forall i &\in \GENERATORS
              \label{model:acopf:pgbound} \\
          & {-}\overline{S_{e}} \leq  \PF_{e} \leq \overline{S_{e}},\;{-}\overline{S_{e}} \leq  \QF_{e} \leq \overline{S_{e}},\;{-}\overline{S_{e}} \leq  \PT_{e} \leq \overline{S_{e}},\;{-}\overline{S_{e}} \leq  \QT_{e} \leq \overline{S_{e}}
              & \forall e &\in \EDGES
              \label{model:acopf:pfbound}
      \end{align}
  \end{subequations}
\end{model}

\clearpage
\subsection{DC-OPF KKT conditions}\label{sec:dckkt}

(a) Primal Feasibility:
These conditions ensure that the primal problem constraints are satisfied:
\begin{align*}
    & \sum_{j \in \GENERATORS_i}\PG_j - \sum_{e \in \mathcal{E}_{i}} \PF_{e} + \sum_{e \in \mathcal{E}^R_{i}} \PF_{e} + \varphi_i
    = \sum_{j \in \LOADS_i}\PD_j + \GS_i, \quad \forall i \in \NODES \\
    & \PF_{e} = -b_e (\VA_i - \VA_j), \quad\quad\quad\quad\quad\quad\quad\quad\enspace\quad\; \forall e = (i, j) \in \EDGES \\
    & \dvamin_e \leq \VA_i - \VA_j \leq \dvamax_e, \quad\quad\quad\quad\quad\quad\enspace\quad\quad \forall e = (i, j) \in \EDGES \\
    & \VA_{\text{ref}} = 0, \quad\quad\quad\quad\quad\quad\quad\enspace\quad\quad\quad\quad\enspace\quad\quad\quad \text{(reference bus)} \\
    & \underline{\PG_i} \leq \PG_i \leq \overline{\PG_i}, \quad\quad\quad\quad\quad\quad\quad\quad\;\,\quad\quad\quad\enspace\;\quad\quad\quad\; \forall i \in \GENERATORS \\
    & -\overline{S_e} \leq \PF_{e} \leq \overline{S_e}, \quad\quad\quad\quad\quad\quad\quad\;\quad\quad\quad\enspace\;\quad\quad\quad \forall e \in \EDGES
\end{align*}

(b) Stationarity:
Stationarity sets the Lagrangian derivatives with respect to primal variables to zero:
\[
\begin{aligned}
    &\frac{\partial \mathcal{L}}{\partial \PG_j} = c_j + \lambdaP_i - \muPgMin_j + \muPgMax_j = 0 \quad \forall j \in \GENERATORS, \\
    &\frac{\partial \mathcal{L}}{\partial \PF_e} = -\lambdaP_i + \lambdaP_j + \lambdaPf_e - \muPfMin_e + \muPfMax_e = 0 \\
    &\quad\quad\quad \forall e = (i, j) \in \EDGES, \\
    &\frac{\partial \mathcal{L}}{\partial \VA_i} = \sum_{e \in \mathcal{E}_i} b_e \lambdaPf_e 
    - \sum_{e \in \mathcal{E}^R_i} b_e \lambdaPf_e + \sum_{e \in \EDGES_i} (-\muAngleDiffMin_e + \muAngleDiffMax_e) = 0 \\
    &\quad\quad\quad \forall i \in \NODES, \\
    &\frac{\partial \mathcal{L}}{\partial \varphi_i} = c_\varphi + \lambdaP_i = 0 \quad \forall i \in \NODES.
\end{aligned}
\]

(c) Dual Feasibility:
Dual variables associated with inequality constraints must remain non-negative:
\begin{align*}
    & \muPgMin_j \geq 0, \quad \muPgMax_j \geq 0, & \forall j \in \GENERATORS \\
    & \muPfMin_e \geq 0, \quad \muPfMax_e \geq 0, & \forall e \in \EDGES \\
    & \muAngleDiffMin_e \geq 0, \quad \muAngleDiffMax_e \geq 0, & \forall e \in \EDGES
\end{align*}

(d) Complementary Slackness:
Complementary slackness ensures that either the inequality constraint is active or its corresponding dual variable is zero:
\begin{align*}
    & \muPgMin_j \perp (\underline{\PG_j} - \PG_j), & \forall j \in \GENERATORS \\
    & \muPgMax_j \perp (\PG_j - \overline{\PG_j}), & \forall j \in \GENERATORS \\
    & \muPfMin_e \perp (-\overline{S_e} - \PF_e), & \forall e \in \EDGES \\
    & \muPfMax_e \perp (\PF_e - \overline{S_e}), & \forall e \in \EDGES \\
    & \muAngleDiffMin_e \perp (\dvamin_e - (\VA_i - \VA_j)), & \forall e = (i, j) \in \EDGES \\
    & \muAngleDiffMax_e \perp ((\VA_i - \VA_j) - \dvamax_e), & \forall e = (i, j) \in \EDGES
\end{align*}

\clearpage
\section{Experiment Setup and Additional Results}
\begin{figure}[!ht]
    \centering
    \includegraphics[width=0.6\linewidth]{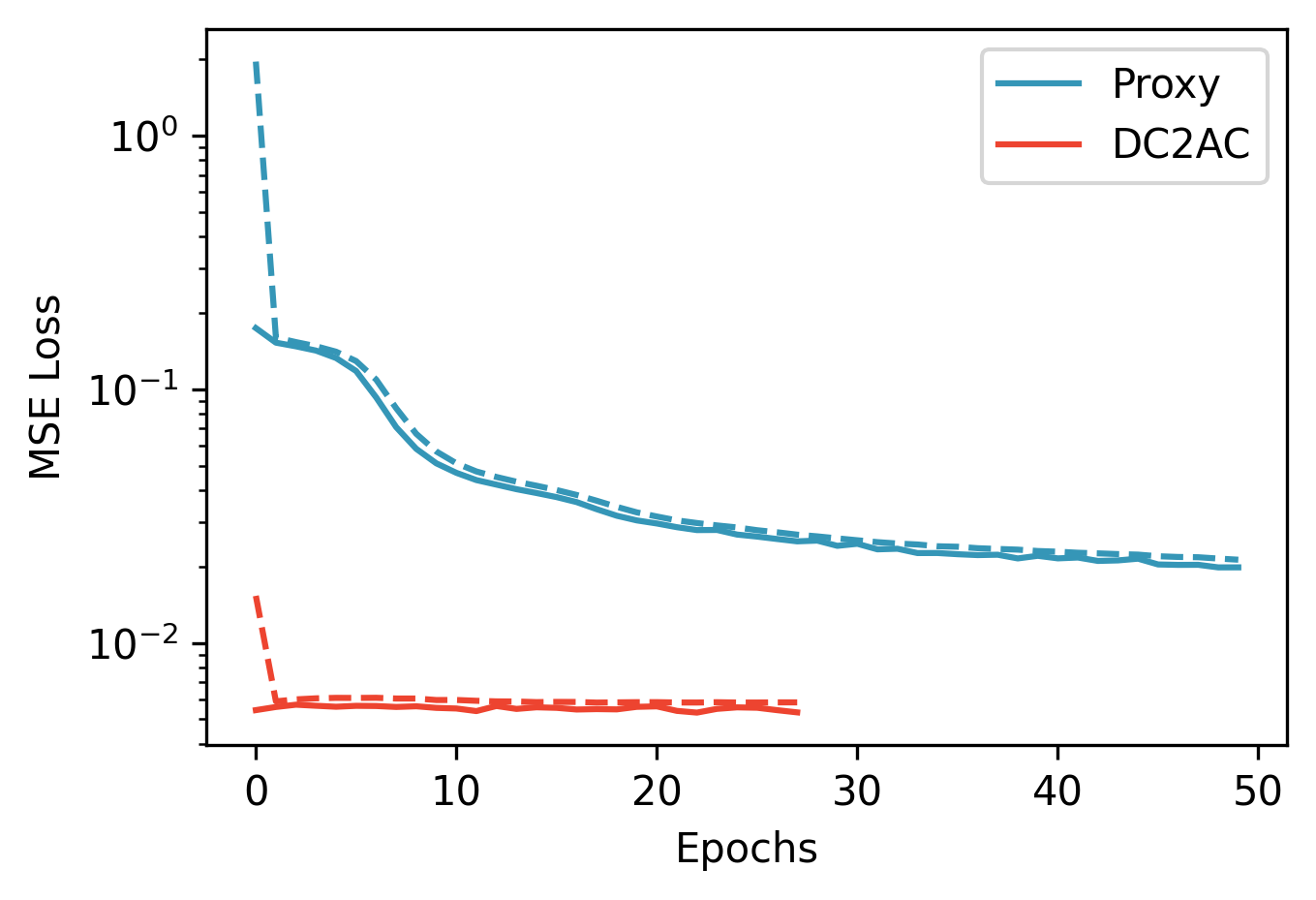}
    \caption{Training (dashed) and validation (solid) loss convergence throughout training for Proxy and DC2AC.}
    \label{fig:conv}
\end{figure}

This section includes details on experiment setup and additional results, namely the convergence plot in Figure \ref{fig:conv} and prediction accuracies in Table \ref{tab:accuracy}.

\subsection{Experiment Setup}\label{sec:expdetail}
The baseline and the proposed method both use a feed-forward neural network architecture with 3 layers each of width 64 and Softplus activations, a constant learning rate of $10^{-4}$, and the Adam \cite{adam} optimizer. Bounds on predictions are enforced using a double-sided Softplus, i.e. $\operatorname{softplus}(x-l)-\operatorname{softplus}(x-u)+l$ to enforce $l\leq x\leq u$. The differentiable optimization layer implementation uses the Julia programming language \cite{julia}; its main dependencies are the JuMP modeling language \cite{jump} and the DiffOpt differentiable optimization framework \cite{diffopt}. The PGLearn.jl library \cite{opfgenerator} is used to build the JuMP model. The learning approaches are implemented in PyTorch \cite{torch} using ML4OPF \cite{opfgenerator}, with DC2AC relying on Julia interface \texttt{juliafunction} \cite{juliafunction} and the LP solver HiGHS \cite{highs}.

\subsubsection{Data Generation Procedure}
This work uses PGLearn.jl \cite{opfgenerator} for data generation, sampling loads using both a global (one per sample) correlation factor as well as local noise (one per sample per load). Both are sampled from Uniform distributions; the local noise range is set to $\pm15\%$ around the each load's reference value and the global range is set to $70\%-110\%$. This ensures the dataset captures a wide range of operating conditions, especially in the challenging high-load regime.
Results are provided for a dataset of 10,000 samples ($80\%$ training, $20\%$ validation) generated around \texttt{89\_pegase} \cite{pegase}, a  benchmark based on the European power grid that has 89 buses (nodes), 35 loads, 12 generators, and 210 transmission lines/transformers (edges).

\begin{table}[ht]
    \centering
    \caption{Prediction accuracies compared to AC-OPF (per-unit).}
    \begin{tabular}{c|ccc}
    \toprule
        Method & $\|\tilde{\mathbf{p}}^\text{g}-\mathbf{p}^\text{g}\|_1$ & $\|\tilde{\mathbf{p}}^\text{f}-\mathbf{p}^\text{f}\|_1$ & $\|\tilde{\VA}-\VA\|_1$ \\
        \midrule
        DC-OPF       & 0.90        & 9.72        & 1.12          \\
        Proxy     & 0.74        & 12.83        & \textbf{{0.95}}          \\
        DC2AC     & \textbf{{0.22}}        & \textbf{{7.57}}        & 1.08          \\
    \end{tabular}
    \label{tab:accuracy}
\end{table}

\end{document}